\documentclass{ifacconf}

\usepackage{graphicx}      
\usepackage{amsmath} 
\usepackage{amssymb}  
\usepackage{amstext,enumerate}
\usepackage{subfigure}
\usepackage{epstopdf}
\usepackage{url}
\usepackage{multirow}
\usepackage{latexsym, color,booktabs}
\usepackage{circuitikz}

\usepackage{natbib}        
\begin{document}
\begin{frontmatter}

\title{Adaptive voltage regulation of an inverter-based power distribution network with a class of droop controllers\thanksref{footnoteinfo}} 

\thanks[footnoteinfo]{This work is supported in part by the Swedish Research Council (grant 2016-0861), the Swedish Energy Agency (project LarGo!), and the Swedish Civil Contingencies Agency (project CERCES).}
\thanks[footnoteinfo]{ © 2020 Michelle S. Chong and Henrik Sandberg. This work has been accepted to IFAC for publication under a Creative Commons Licence CC-BY-NC-ND}

\author[tue]{Michelle S. Chong} 
\author[kth]{Henrik Sandberg}

\address[tue]{Department of Mechanical Engineering, Eindhoven University of Technology, the Netherlands. {\tt\small m.s.t.chong@tue.nl }}
\address[kth]{Division of Decision and Control Systems, KTH Royal Institute of Technology, Sweden.
        {\tt\small hsan@kth.se}}

\begin{abstract}                
 The voltage received by each customer connected to a power distribution line with local controllers (inverters) is regulated to be within a desired margin through a class of slope-restricted controllers, known conventionally as \emph{droop} controllers. We adapt the design of the droop controllers according to the known bounds of the net power consumption of each customer in each observation time window. A sufficient condition for voltage regulation is provided for each time window, which guides the design of the droop controllers, depending on the properties of the distribution line (line impedances) and the upper bound of all the customers' power consumption during each time window. The resulting adaptive scheme is verified on a benchmark model of a European low-voltage network by the CIGRE task force.   
\end{abstract}


\end{frontmatter}

\section{Introduction}
The proliferation of harnessing renewable energy sources such as solar and wind at the customer end has been motivated by meeting environmental sustainability goals in lieu of being financially viable. This is actuated via connecting alternating current/direct current (AC/DC) inverters which are co-located with the customer to the power distribution network, see \cite{lopes2007integrating}. As a result, the customer departs from the traditional role of being merely a consumer, as the customer can now have a dual role of also being a producer. Consequently, fluctuations outside of the safety margin in the voltage received by the customer can occur, especially when there are many customers. Therefore, this calls for a systematic design of the local inverters (controllers) to guarantee that the voltages at the customer end are within a safety margin specified by the operator. 

This paper concentrates on regulating the voltage level of each customer connected to a power distribution network in a line configuration, which is achieved by adapting the controller design based on each customer's projected power consumption. In this setup, power is delivered to each customer with the substation at the head of the line, and the customers connect to line sequentially, as shown in Figure \ref{fig:ems}. A nominal voltage $\bar{v}$ is communicated by the substation to each inverter, which injects reactive power by employing the so-called droop controllers, which has been used for voltage regulation in \cite{andren2015stability,turitsyn2010distributed,simpson2017voltage,farivar2013equilibrium,jahangiri2013distributed,han2017containment,zhu2016fast,liu2008distribution,schiffer2014conditions,jafarian2018interconnection}, to name a few. 

In our previous papers, \cite{chong2019voltage} and \cite{chong2019local}, we showed that voltage regulation can be guaranteed with a class of slope-restricted droop controllers, which differs from the conventional droop controllers used in the literature listed earlier, in two ways: (i) the input to the controller is the difference of the \emph{squared} received and nominal voltages, and (ii) a class of droop functions satisfying a slope restriction property is identified, which expands beyond the piecewise continuous saturation function used in the literature. Further, a sufficient condition for voltage regulation uses only the upper bound of all the customers' net power consumption, which has the additional benefit of preserving privacy. A drawback of our previous results is a sufficient condition that is overly conservative, resulting in a higher than needed generation of reactive power by the controllers. In this paper, we consider finite observation time windows, where the upper bounds of the customers' power consumption is known beforehand. This allows us to adapt the design of our class of droop controllers, while still guaranteeing that all the customer's voltage are within a desired margin of the nominal voltage. This adaptive scheme is tested on a benchmark model of the residential feeder of a European low-voltage distribution line, and shown in simulations that less reactive power needs to be generated by the controllers. Additionally, more active power can be injected into the line, which increases the revenue of each customer. 

The remainder of the paper is structured as follows. We first introduce the notations used, followed by Section \ref{sec:model} where we present a model of a power distribution line with inverters. The voltage regulation objective is then stated in Section \ref{sec:problem}, which is achieved via local injection of reactive power into the distribution network by the inverters. In Section \ref{sec:main}, we present a scheme that adapts each local controllers' droop characteristics according to the a priori known bounds on the net active power and consumed reactive power during an observation window. We then verify this adaptive scheme on a benchmark model in Section \ref{sec:sim} and show in simulations that the adaptive scheme performs better than the non-adaptive scheme according to two performance metrics. Section \ref{sec:conc} concludes the paper with discussions for future work.

\subsubsection{Notation} 
 	 Let $\mathbb{R}=(-\infty,\infty)$, $\mathbb{R}_{\geq 0}=[0,\infty)$, $\mathbb{R}_{>0}=(0,\infty)$, $\mathbb{N}=\{0,1,\dots\}$.
	 Let the set of complex numbers be denoted by $\mathbb{C}$.
	 We denote the set of integers $\{i,i+1,i+2,\dots,i+k\}$ as $\mathbb{N}_{[i,i+k]}$.
 	 Let $(u,v)$ where $u\in\mathbb{R}^{n_u}$ and $v\in\mathbb{R}^{n_v}$ denote the vector $(u^T,v^T)^{T}$.
 	 The identity matrix of dimension $n$ is denoted by $\mathbb{I}_{n}$ and a matrix of dimension $m$ by $n$ with all elements $1$ is denoted by $\mathbf{{1}}_{m\times n}$.
	 A diagonal matrix with elements $d_i$, $i\in\mathbb{N}_{[1,n]}$ is denoted by $\textrm{diag}(d_1,d_2,\dots,d_n)$.
 	 Given a symmetric matrix $P$, its maximum (minimum) eigenvalue is denoted by $\lambda_{\max}(P)$ $(\lambda_{\min}(P))$.
	 The infinity norm of a vector $x \in \mathbb{R}^{n}$, is denoted $|x|:= \underset{i\in\mathbb{N}_{[1,n]}}{\max} \left| x_i \right|$ and for a matrix $A\in\mathbb{R}^{n\times n}$, $|A|:= \underset{i\in\mathbb{N}_{[1,n]}}{\max} \underset{j\in\mathbb{N}_{[1,n]}}{\sum}|a_{ij}|$, where $a_{ij}$ is the row $i$-th and column $j$-th element of matrix $A$.

\section{Model of a power distribution line with inverters} \label{sec:model}
\begin{figure}[h!]
	\begin{center}
		\includegraphics[scale=0.85]{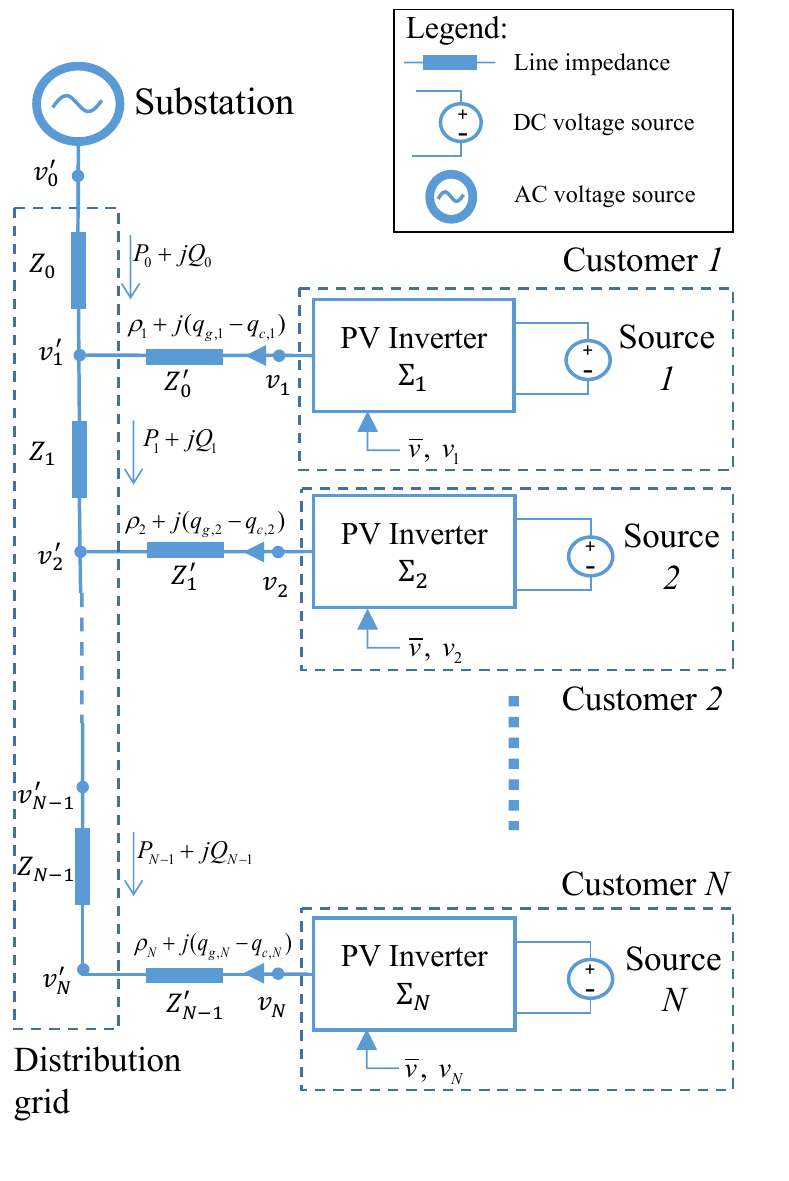}
	\end{center} \vspace{-2em}
	\caption{Infrastructure of the low-voltage distribution grid. } \label{fig:ems}
\end{figure}

We consider a model of a distribution grid infrastructure with $N$ customers feeding into a low-voltage (LV) grid in a line configuration, as shown in Figure \ref{fig:ems}. Each customer $i\in\mathbb{N}_{[1,N]}$ has a renewable energy source such as a photovoltaic cell, modelled as a controllable voltage source \citep{simpson2017voltage}. Each renewable energy source is equipped with an inverter\footnote{This is enforced by curtailing the generated active power $\rho_{g,i}$ based on the generated reactive power $q_{g,i}$ and the power rating of the inverter, known as the apparent power $\bar{s}_{i}$, according to $\rho_{g,i} \leq \sqrt{\bar{s}_{i}^2-q_{g,i}^2}$} which is able to generate reactive power $q_{g,i}$ while generating near maximum active power $\rho_{g,i}$. 

For each customer $i\in\mathbb{N}_{[1,N]}$, the received voltage level is $v_i$ and the voltage level at the point of connection with the distribution line is $v'_i$, with a corresponding line impedance $Z'_i=R'_{i} + jX'_{i}$ in between customer $i$ and the connection point on the distribution line, where $R'_{i}\in\mathbb{R}_{\geq 0}$ is the resistance and $X'_{i}\in\mathbb{R}_{\geq 0}$ is the reactance. In between each connection point, the corresponding line impedance is $Z_{i}=R_{i} + jX_{i}$, where $R_{i}\in\mathbb{R}_{\geq 0}$ is the resistance and $X_{i}\in\mathbb{R}_{\geq 0}$ is the reactance. Each customer has a load which can consume reactive $q_{c,i}$ and active powers $\rho_{c,i}$, independently of the generated reactive $q_{g,i}$ and active powers $\rho_{g,i}$.

Under the assumption that the power losses in the distribution lines are negligible, we use the linearized DistFlow model in \cite{baran1989network} for describing the power flow in a line configuration as the nonlinear terms capture the power losses in the distribution lines. Additionally, we also model the point-to-point voltages of the connection point at the distribution line $v'_i$ and the customer's end $v_i$ in the last equation of \eqref{eq:line_rel}. The relationship between the power flow and point-to-point voltage of key nodes of interest $i\in \mathbb{N}_{[0,N]}$, in Figure \ref{fig:ems} is
\begin{equation} \label{eq:line_rel}
	\begin{array}{lll}
			P_{i+1} & = & P_{i} + \rho_{i+1}, \\
			Q_{i+1} & = & Q_{i} + q_{i+1}, \\
			{v'_{i+1}}^2 & = & {v'_{i}}^2 - 2\beta_{i}(P_i,Q_i), \\
			{v'_i}^2 & = & {{v}_{i}}^2 - 2 {\beta}^{'}_{i-1}(\rho_i,q_i),
	\end{array}
\end{equation}
where $P_i$ and $Q_i$ are the respective total active and reactive powers flowing from customer $i$ to customer $i+1$; $\rho_i:=\rho_{g,i}-\rho_{c,i}$ and $q_i:=q_{g,i} - q_{c,i}$ are the net injection of the respective active and reactive power into the distribution line from customer $i$; $\beta_{i}(r,s):=R_i r + X_i s$ and $\beta'_{i}(r,s):=R'_i r + X'_i s$ with $\beta'_{-1}(r,s)=0$ for all $r,s\in\mathbb{R}$.  

Our modelling assumptions are stated in Assumption \ref{assum:model} below.
\begin{assum}[Modelling assumptions] \label{assum:model}
	For each \\ customer $i\in \mathbb{N}_{[1,N]}$:
	\begin{enumerate}[(i)]
		\item The reactive power $q_{c,i}$ consumed by each customer $i$ is bounded and the bound $\Delta_c$ is known, i.e. $|q_{c,i}(t)|\leq \Delta_c$, for all $t\geq 0$.
		\item The net injected active power $\rho_{i}$ is bounded and the bound $\Delta_\rho$ is known, i.e. i.e. $|\rho_{i}(t)|\leq \Delta_\rho$, for all $t\geq 0$.
		\item The node voltages at the receiving end $v_i$ for each customer $i$ is locally measurable.
		\item The nominal voltage reference $\bar{v}$ is communicated to each customer $i$. Standard values for a single-phase and three-phase systems are $230$ V and $400$ V, respectively.
	\end{enumerate} 
\end{assum}

\section{Problem formulation} \label{sec:problem}
We aim to regulate the customers' received voltage levels $v_i$ for $i\in\mathbb{N}_{[1,N]}$, such that it lies within a safe operating range of the reference voltage communicated to each customer, $\bar{v}\in\mathbb{R}$, i.e. for a given $\delta\geq 0$, 
\begin{equation}
	\bar{v}-\delta \leq v_{i}(t) \leq \bar{v}+\delta, \qquad \forall t\geq 0. \label{eq:ori_voltage_bounds} 
\end{equation}
This is achieved by only managing the flow of the reactive powers $Q_{i}$ in the distribution grid. To this end, controllers $\Sigma_i$ (actuated by the inverters) are designed for each customer, which is also known in the literature such as in \cite{andren2015stability}, as droop control. In \cite{chong2019local}, we showed that a class of slope-restricted droop controllers can be designed to achieve the voltage regulation objective \eqref{eq:ori_voltage_bounds}, for the power distribution model considered in this paper. 

In this paper, we consider the case where we update the droop controllers $\Sigma_i$ at discrete time instants $t_k$, $k\in\mathbb{N}$, due to new information about the bounds on the power consumption ($\Delta_c$ and $\Delta_\rho$ respectively. Cf. Assumption \ref{assum:model}). Without loss of generality, we consider $t_0 = 0$ and assume that the observation time window $[t_k,t_{k+1})$ is sufficiently large such that no Zeno phenomenon may occur. Further, the bounds are known a priori for each observation window $[t_k,t_{k+1})$. For example, power consumption in a commercial district peaks during business hours and declines after hours. The motivation behind this setup is to adapt our controllers to the known bounds on the power consumption for certain time epochs. We denote these bounds on the consumed reactive and active powers during the observation window $[t_k,t_{k+1})$ as $\Delta_{c,k}$ and $\Delta_{\rho,k}$, respectively. We then update the droop function $K_{i,k}$ accordingly, in the droop controllers $\Sigma_i$ defined as follows
\begin{equation} \label{eq:droop}
\Sigma_i: \; \dot{q}_{g,i} = -\frac{1}{\tau_i} q_{g,i} + \frac{1}{\tau_i} K_{i,k}(\bar{v}^2-v_{i}^2), \; \forall t\in[t_k,t_{k+1}), 
\end{equation} 
where $\tau_i\in\mathbb{R}_{> 0}$ is the time-constant of the inverter's response and the droop function $K_{i,k}(w)$ is a static mapping from the difference of the squared voltages $w$ to the set-point for the reactive power, updated at $t_{k}$. The droop function $K_{i,k}$ satisfies the slope-restriction and saturation property \eqref{eq:sector_bound} as assumed in \cite{chong2019local}, which we state in the following assumption.

\begin{assum}\label{assum:sector_bound}
	During each observation window $t\in[t_k,t_{k+1})$, $k\in\mathbb{N}$, the droop function $K_{i,k}$, $i\in\mathbb{N}_{[1,N]}$ in \eqref{eq:droop} satisfies
	\begin{equation} \label{eq:sector_bound}
			0 \leq \frac{K_{i,k}(v)-K_{i,k}(w)}{v-w} \leq d_{i,k}, \qquad \forall v,w\in\mathbb{R}, \; v\neq w,
	\end{equation}
	and $|K_{i,k}(v)|\leq \bar{K}_{i,k}$, for all $v\in\mathbb{R}$, where $d_{i,k} \geq 0$ and $\bar{K}_{i,k}\geq 0$. \hfill $\Box$
\end{assum}

The parameter $d_{i,k}$ is a design parameter chosen according to the condition in Section \ref{sec:main}. 

At each update time $t_k$, we reset the state of all the controllers \eqref{eq:droop} to zero\footnote{This choice was made to simplify the re-design of the droop function $K_{i,k}$ according to condition \eqref{eq:cond_volt_reg} in Theorem \ref{thm:voltage_reg}.} as follows
\begin{equation}
	q_{g,i}(t_{k}^{+}) = 0. \label{eq:init_condition}
\end{equation}

Under this setup, we aim to meet the voltage regulation objective \eqref{eq:ori_voltage_bounds}.

\section{Adaptive scheme for voltage regulation} \label{sec:main}
To assist in formulating our main result, we rewrite the distribution model \eqref{eq:line_rel} and droop controllers \eqref{eq:droop} in new coordinates, with $y_i=\bar{v}^2 - v_i^2$, $\eta_i:=v_{i-1}^2 - v_{i}^2$ for $i\in\mathbb{N}_{[1,N]}$ and $y_0:=\bar{v}^2 - {v'}_{0}^{2}$. The resulting system for $i\in\mathbb{N}_{[1,N]}$ is 
\begin{equation}
	 \begin{array}{rll}
			\dot{q}_{g,i} & = & -\frac{1}{\tau_i} q_{g,i} + \frac{1}{\tau_i}K_{i,k} (y_i), \qquad t\in[t_{k},t_{k+1}), \\
			q_{g,i}(t_k^{+}) & = & 0, \\
			y_i & = & \eta_i + y_{i-1}, \\
			\eta_{i} & = & v_{i-1}^2 -{v'}_{i}^2 + {v'}_{i}^2 - v_i^2, \\
			&= & 2 \beta_{i-1}(P_{i-1},Q_{i-1}) + 2 \beta'_{i-2}(\rho_{i-1},q_{i-1}) \\
			&&- 2 \beta'_{i-1}(\rho_i,q_i), 
	\end{array}  \label{eq:model2}
\end{equation}
where we recall that $\beta_{i}(r,s)=R_{i}r + X_{i}s$, $\beta'_{i}(r,s)=R'_i r + X'_i s$ and $\beta'_{-1}(r,s)=0$ for all $r,s\in\mathbb{R}$.

Our voltage regulation objective \eqref{eq:ori_voltage_bounds} can then be restated in the new coordinates as:
\begin{equation}
	\left|\bar{v}^2-v_i(t)^2\right| = \left|y_{i}(t) \right| \leq \epsilon, \qquad \forall i\in\mathbb{N}_{[1,N]}, \; \forall t\geq 0, \label{eq:restate_voltage_bounds} 
\end{equation}
where $\epsilon:=-\delta^2 + 2\bar{v} \delta$. Note that the satisfaction of \eqref{eq:restate_voltage_bounds} implies that objective \eqref{eq:ori_voltage_bounds} holds. Henceforth, we will state our main result such that that the equivalent objective \eqref{eq:restate_voltage_bounds} is satisfied.

Before doing so, we write the system \eqref{eq:model2} in the compact form below, by applying \cite[Proposition 4]{chong2019local}.
Let $q_{g}:= \left(q_{g,1},q_{g,2},\dots,q_{g,N} \right)$, $q_{c}:= \left(q_{c,1},q_{c,2},\dots,q_{c,N} \right)$, $\rho:=\left(\rho_1,\rho_2,\dots,\rho_{N}\right)$ and $y:=\left(y_1,y_2,\dots,y_N\right)$. The model \eqref{eq:model2} in compact form is written as follows
\begin{equation}
	\begin{array}{rll}
		\dot{q}_{g} & = & A(\tau) q_{g} - A(\tau)K_{k}(y), \qquad t\in[t_{k},t_{k+1}),\\
		q_{g}(t_{k}^{+}) & = & 0, \\
		y & = & Hq_{g} + \phi(\rho,q_c) + y_0 \mathbf{1}_{N\times 1},
	\end{array} \label{eq:model_compact}
\end{equation}  
where 
\begin{itemize}
		\item $A(\tau) = \textrm{diag}\left( -\frac{1}{\tau_1},-\frac{1}{\tau_2},\dots,-\frac{1}{\tau_N}  \right)$,
		\item $K_{k}(y) = \left( K_{1,k}(y_1), K_{2,k}(y_2), \dots, K_{N,k}(y_N) \right)$,
		\item 
		\begin{equation*}
			\begin{aligned}
				H =  -2 & \left(\begin{array}{ccccc} X_0 &  X_0 & \dots  & X_0 \\
			 								  \star & \left(X_0+X_1\right) & \dots  & \left(X_0+X_1 \right) \\
											  \vdots & \ddots & \ddots & \vdots \\
											  \star & \dots & \star & \underset{i\in\mathbb{N}_{[0,N-1]}}{\sum} X_i \end{array}\right)\\ 
											  &-2  \textrm{diag}\left(X'_0,\dots,X'_{N-1} \right),
			\end{aligned}
		\end{equation*}				
		where $\star$ denotes a block component of a symmetric matrix,
											  
		\item $\phi(\rho,q_c)= \left(\phi_{0}, \dots,\!\!\! \underset{i\in\mathbb{N}_{[0,j-1]}}{\sum} \phi_{i}+\!\!\! \underset{i\in\mathbb{N}_{[0,j-2]}}{\sum} 2\beta'_i(\rho_{i+1},q_{c,i+1}) \right.$ $\left.  ,\dots, \underset{i\in\mathbb{N}_{[0,N-1]}}{\sum} \phi_{i} +  \underset{i\in\mathbb{N}_{[0,N-2]}}{\sum} 2\beta'_i(\rho_{i+1},q_{c,i+1}) \right)$, 
		\\ where 
		\begin{equation*}
			\begin{aligned}
				\phi_{i}(\rho,q_c)&:=  2 X_i \left( \underset{j\in\mathbb{N}_{[i+1,N]}}{\sum} q_{c,j} \right)-2 R_i \left( \underset{j\in\mathbb{N}_{[i+1,N]}}{\sum} \rho_j \right)\\ &- 2 \beta'_{i}(\rho_{i+1},q_{c,i+1}).
			\end{aligned}
		\end{equation*}
\end{itemize}

Let $\{t_{k}\}_{k\in\mathbb{N}}$ where $t_{k+1}-t_{k} \geq T >0$, $\{\Delta_{c,k}\}_{k\in\mathbb{N}}$ and $\{\Delta_{\rho,k}\}_{k\in\mathbb{N}}$ be a sequence of update time instants, a sequence of a priori known bounds on the consumed reactive power and net injected active power, respectively. We are now ready to state our main result, where the proof can be found in the Appendix. 

\begin{thm} \label{thm:voltage_reg}
	Consider the distribution model \eqref{eq:line_rel} and droop controllers \eqref{eq:droop} with $N$ customers under Assumption \ref{assum:model} and \ref{assum:sector_bound}.  For any $\epsilon\geq 0$, if there exist $\Delta_{\phi,k}=\Delta_{\phi,k}(\Delta_{\rho,k},\Delta_{c,k})>0$ and $\epsilon_{y}>0$ for all $k\in\mathbb{N}$, such that
	\begin{equation}
		\epsilon \geq \Delta_{\phi,k} + \epsilon_{y},
	\end{equation}
	then there exist $d_{k}\geq 0$, $\bar{K}_{k} \geq 0$, $\tau_{\max} > 0$ and $\tau_{\min}>0$ satisfying 
	\begin{equation}\label{eq:cond_volt_reg}
 				  d_{k} \leq \frac{\left(\epsilon-\Delta_{\phi,k} - \epsilon_y \right) \exp(-t_{k}/\tau_{\max})}{\frac{\tau_{\max}}{\tau_{\min}}  \left|H \right| \left( \Delta_{\phi,k} + \epsilon_y \right)  +  \left(\frac{\tau_{\max}}{\tau_{\min}}\right)^2  \left|H\right|^2  \bar{K}_k},
	\end{equation}
	where 
		\begin{itemize}
			\item $d_k:=\underset{i\in\mathbb{N}_{[1,N]}}{\max} d_{i,k}$, 
			\item $\tau_{\max}:= \underset{i\in\mathbb{N}_{[1,N]}}{\max}\tau_i$ and $\tau_{\min}:= \underset{i\in\mathbb{N}_{[1,N]}}{\min}\tau_i$,
			\item $\bar{K}_k:=\underset{i\in\mathbb{N}_{[1,N]}}{\max} \bar{K}_{i,k}$, 
			\item \begin{equation*}
					\begin{aligned}
						\Delta_{\phi,k}&:={N(N+1)}\left( \bar{R} \Delta_{\rho,k} + \bar{X} \Delta_{c,k} \right)\\
						& \qquad + 2(2N-1) \left( \bar{R'} \Delta_{\rho,k} +\bar{X'}\Delta_{c,k} \right),
					\end{aligned}
					\end{equation*}
			where
			\begin{itemize}
				\item $\bar{R}:=\underset{i\in\mathbb{N}_{[1,N-1]}}{\max} R_i$, $\bar{R'}:=\underset{i\in\mathbb{N}_{[0,N-2]}}{\max} R'_i$,
				\item $\bar{X} := \underset{i\in\mathbb{N}_{[1,N-1]}}{\max} X_i$, $\bar{X'} := \underset{i\in\mathbb{N}_{[0,N-2]}}{\max} X'_i$,
			\end{itemize}
			\item $|\bar{v}^2-v'_{0}(t)^2|\leq \epsilon_y$, for all $t\geq 0$,
	\end{itemize}		  
	such that the received voltage $v_i$ of each customer $i$ meets the voltage regulation objective \eqref{eq:restate_voltage_bounds}. \hfill$\Box$
\end{thm}

The sufficient condition \eqref{eq:cond_volt_reg} in Theorem \ref{thm:voltage_reg} provides a guideline for adapting the droop controller \eqref{eq:droop} according to the a priori known bounds $\Delta_{c,k}$, $\Delta_{\rho,k}$ on the power consumption during each observation window $[t_{k},t_{k+1})$, $k\in\mathbb{N}$. The design parameter $d_{k}$ in \eqref{eq:cond_volt_reg} guides the design of the droop function $K_{i,k}$ in each of the droop controllers \eqref{eq:droop}. For instance, a piecewise saturation function commonly considered in the literature such as \cite{andren2015stability} takes the form below
	\begin{equation} \label{eq:droop_K}
	\begin{aligned}
			& K_{i,k}(w):= \\ &  \left\{\begin{array}{ll} -\bar{Q}_{i,k}, & w \leq w_{\min,i,k}, \\  
			  -\left(1-  \frac{w-w_{\min,i}}{w_{m,i,k}-w_{\min,i,k}}\right) \bar{Q}_{i,k},  & w\in (w_{\min,i,k},w_{m,i,k}], \\
			 0, & w\in (w_{m,i,k},w_{n,i,k}], \\
		  \left(\frac{w-w_{n,i,k}}{w_{\max,i,k}-w_{n,i,k}}\right) \bar{Q}_{i,k},  &  w\in (w_{n,i,k},w_{\max,i,k}], \\
			 \bar{Q}_{i,k}, & w > w_{\max,i,k},  \end{array} \right. 
	\end{aligned}
	\end{equation}
where $w_{\min,i,k} \leq w_{m,i,k} \leq 0 \leq w_{n,i,k} \leq w_{\max,i,k}$ are design parameters, $\bar{Q}_{i,k}\in\mathbb{R}_{\geq 0}$ is the saturation limit of the $i$-th inverter satisfying $\bar{Q}_{i,k}=\bar{s}_i$, where $\bar{s}_i\in\mathbb{R}$ is the maximum apparent power of the $i$-th inverter. This droop function \eqref{eq:droop_K} satisfies Assumption \ref{assum:sector_bound} with $\bar{K}_{i,k}=\bar{Q}_{i,k}$ and
	\begin{equation} \label{eq:di_design}
		d_{i,k}:=\min \left\{\begin{array}{cc} \frac{\bar{Q}_{i,k}}{w_{\max,i,k} - w_{n,i,k}}, \frac{\bar{Q}_{i,k}}{w_{m,i,k} - w_{\min,i,k}} \end{array} \right\}.
	\end{equation}	

\section{Case study on a benchmark model} \label{sec:sim}
\begin{figure}[h!]
    \centering
    \includegraphics[width=\columnwidth]{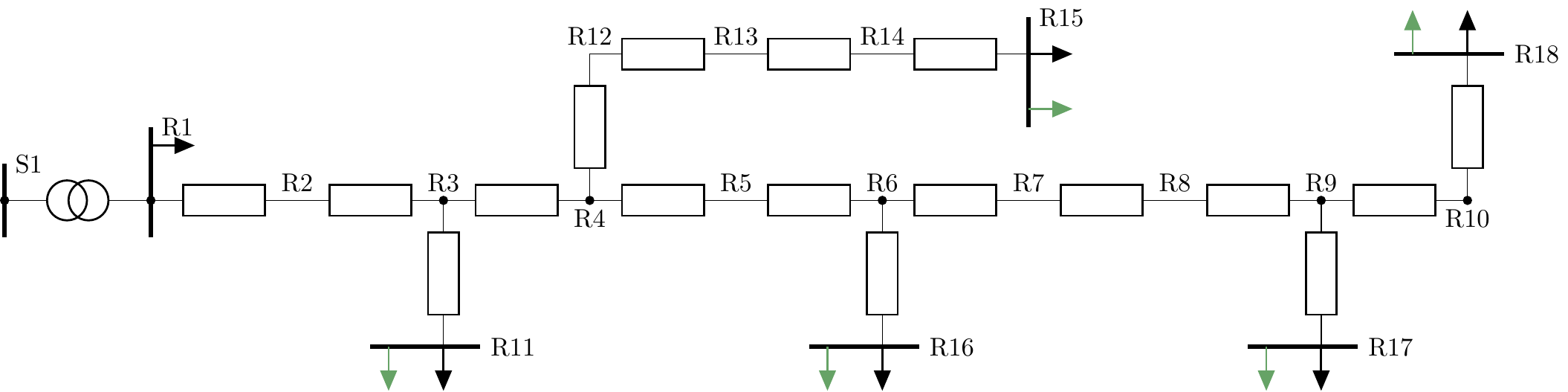}
    \caption{Residential feeder of the European low voltage CIGRE benchmark grid. Arrows in \textbf{{black}} represent loads and the arrows in \textbf{\color{green}{green}} represent inverters.}
    \label{fig:cigre}
\end{figure}

We validate our results on a benchmark model of $N=5$ customers in the residential feeder of the European low voltage CIGRE distribution grid \cite[Fig. 7.7]{strunz2009benchmark}, shown in Figure \ref{fig:cigre}. The topology of the benchmark model can be mapped to our model shown in Figure \ref{fig:ems} according to Table \ref{tab:map_nodes}. The model parameters as found in \cite[Table 7.26]{strunz2009benchmark} are summarised in Table \ref{tab:grid_parameters}. The reference voltage $\bar{v}$ communicated to each customer is $\bar{v}=230$ V and the nominal voltage at the substation is $v'_0(t)=\bar{v}+ 5\sin(t)$ V to model harmonic perturbations.

\begin{table}[h!]
	\centering
	\caption{Mapping of nodes from the benchmark topology to our topology}
	\label{tab:map_nodes}
	\begin{tabular}{|c|c|}
		\hline 
		\multicolumn{2}{|c|}{Node label} \\
		\hline 
		Benchmark topology & \\ \cite[Fig. 7.7]{strunz2009benchmark} & Our topology in Figure \ref{fig:ems} \\ in Figure \ref{fig:cigre}   & \\ 
		\hline 
		$R1$ & $v'_0$ \\ 
		\hline 
		$R3$  & $v'_1$ \\ 
		\hline 
		$R11$  & $v_1$ \\ 
		\hline 
		$R4$  & $v'_2$ \\ 
		\hline 
		$R15$  & $v_2$ \\ 
		\hline
		$R6$  & $v'_3$ \\ 
		\hline 
		$R16$  & $v_3$ \\ 
		\hline 
		$R9$  & $v'_4$ \\
		\hline
		$R17$  & $v_4$ \\
		\hline
		$R10$  & $v'_5$ \\
		\hline
		$R18$  & $v_5$ \\
		\hline
	\end{tabular} 
\end{table}

\begin{table}[h!]
	\centering
	\caption{Parameters as found in \cite[Table 7.26]{strunz2009benchmark}}
	\label{tab:grid_parameters}
	\begin{tabular}{|l|c|c|c|c|c|}
		\hline 
		$i$ & 1 & 2 & 3 & 4 & 5 \\ 
		\hline 
		$R_{i-1}\,[\Omega]$ & 0.00343 & 0.00172 & 0.00343 & 0.00515 & 0.00172 \\ 
		\hline 
		$X_{i-1}\,[\Omega]$ & 0.04711 & 0.02356 & 0.04711 & 0.07067 & 0.02356 \\ 
		\hline 
		$R'_{i-1}\,[\Omega]$ & 0.00147 & 0.00662 & 0.00147 & 0.00147 & 0.00147 \\ 
		\hline 
		$X'_{i-1}\,[\Omega]$ & 0.02157 & 0.09707 & 0.02157 & 0.02157 & 0.02157 \\ 
 		\hline 
		Apparant  & & & & & \\
		power of & & & & & \\
		each  & 4200 & 6500 & 4700 & 5300 & 3600 \\
		inverter & & & & & \\
		$\bar{s}_{i}$, [VA] & & & & & \\
		\hline
	\end{tabular} 
\end{table}

\begin{figure}[h!]
    \centering
    \includegraphics[scale=0.8]{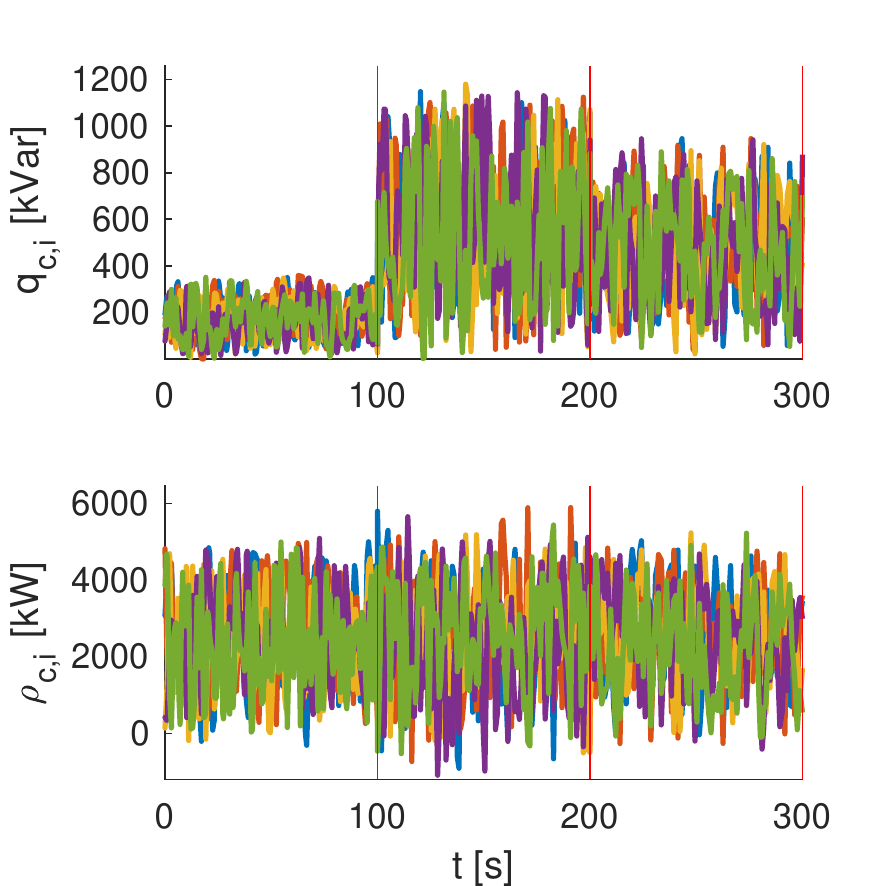}
	\caption{Reactive $q_{c,i}$ and active $\rho_{c,i}$ powers consumed by each customer $i$.} \label{fig:data}
\end{figure}

\begin{figure}[h!] 
	\begin{center}
		\includegraphics[scale=0.8]{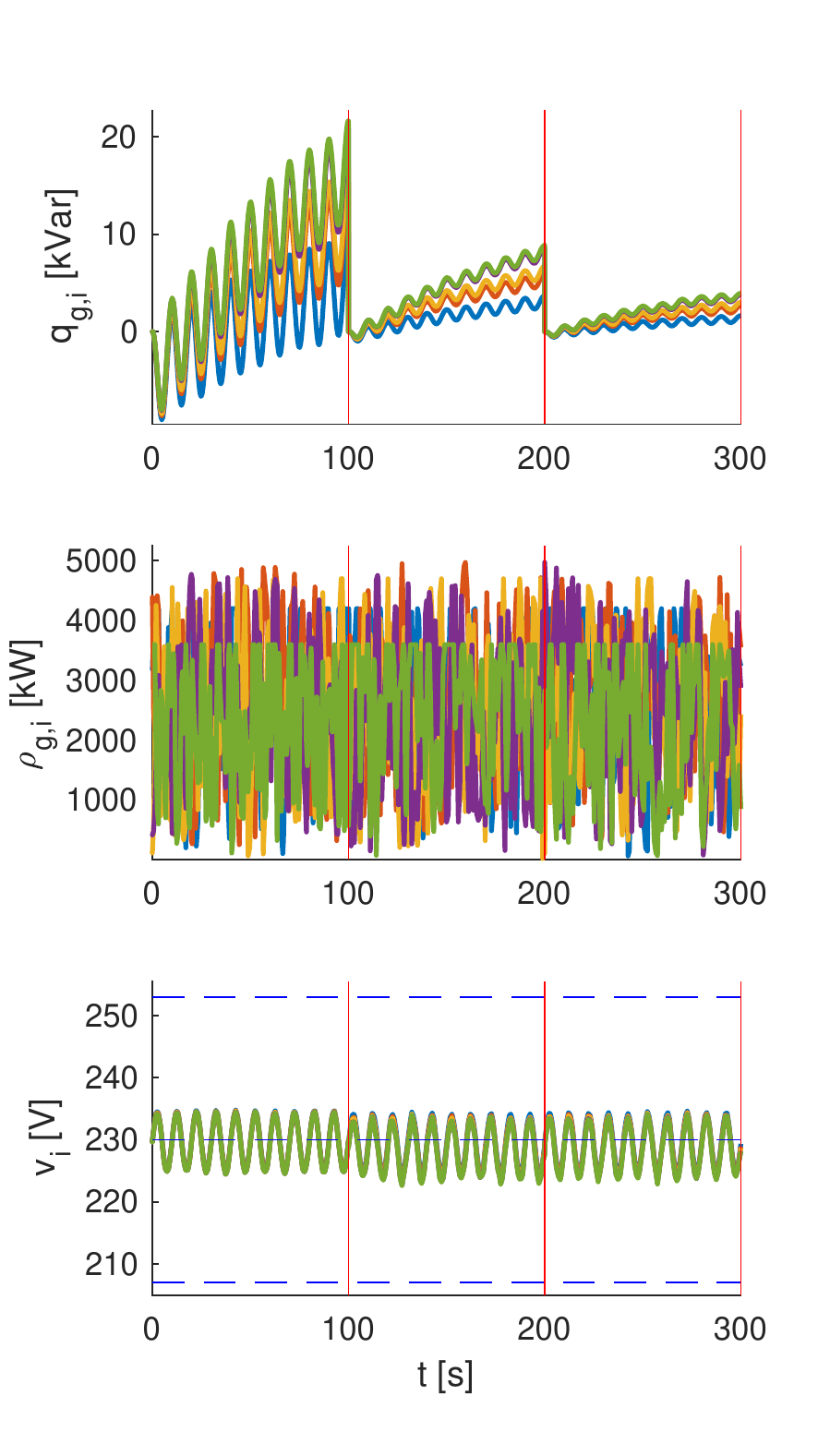} \vspace{-2em}
		\caption{Adaptive droop controllers: Vertical {\color{red}red} lines indicate the observation instants $t_k$, $k\in\mathbb{N}_{[0,2]}$. Legend: Customer $1$ (blue), $2$ (red), $3$ (orange), $4$ (purple), $5$ (green).  (Top) Solid lines indicate the output of the droop controllers $q_{g,i}$. (Middle) The active power $\rho_{g,i}$ injected into the distribution line by each customer $i$. (Bottom) Solid lines indicate the voltages of each customer $v_i$ are within $10$ \% of the nominal voltage $\bar{v}$ for all $t\geq 0$. Dashed lines indicate $(1\pm p)\bar{v}$, where $p=0.1$.  }\label{fig:adapt_time}
	\end{center}
\end{figure}

Our simulation study considers the observation instants $t_0=0$ s, $t_1=100$ s, $t_2 = 200$ s and $t_3 = 300$ s, where the a priori known upper bounds on the net active power and the consumed reactive power has changed in between the observation instants as shown in Table \ref{tab:sim} and as seen in Figure \ref{fig:data}. All the inverters \eqref{eq:droop} have a time constant of $\tau_i=100$, for $i\in\mathbb{N}_{[1,N]}$.  By applying Theorem \ref{thm:voltage_reg}, we choose the design parameter $d_k$ according to \eqref{eq:cond_volt_reg} which then guides the design of the droop functions $K_{i,k}$ in \eqref{eq:droop_K} by choice of the droop characteristics $w_{\min,i,k}$, $w_{m,i,k}$, $w_{n,i,k}$ and $w_{\max,i,k}$. Here, the saturation value $\bar{K}_{i,k}$ is set to the apparent power $\bar{s}_{i}$ of each inverter, $\bar{K}_{i,k}=\bar{s}_{i}$, for all $k\in\mathbb{N}$ (c.f. Table \ref{tab:grid_parameters} for values of $\bar{s}_{i}$). Hence, by the physical capabilities of the inverter, the injected active power is limited to $\rho_{g,i} \leq \sqrt{\bar{s}_{i}^2-q_{g,i}^2}$. Figure \ref{fig:adapt_time} shows that the resulting adaptive scheme achieves the voltage regulation objective \eqref{eq:ori_voltage_bounds}.

\begin{table}[h!]
	\centering
	\caption{Simulation scenario and the application of \eqref{eq:cond_volt_reg} in Theorem \ref{thm:voltage_reg}}
	\label{tab:sim}
	\begin{tabular}{|l|c|c|c|}
		\hline 
		Observation window & & &  \\
		$[t_k,t_{k+1})$ & $k=0$ & $k=1$ & $k=2$ \\ 
		\hline 
		Upper bound on the  & & &  \\
		net injected active  &  &  &  \\
		power for all & & & \\
		customers, $\Delta_{\rho,k}$ & $460$ & $1300$ & $750$ \\ 
		\hline 
		Upper bound on the  & & &  \\
		consumed reactive &  &  &  \\
		power for all & & & \\
		customers, $\Delta_{c,k}$ & $360$ & $1200$ & $960$ \\
		\hline 
		$d_k$ obtained from & & & \\ 
		\eqref{eq:cond_volt_reg} in Theorem \ref{thm:voltage_reg} & $0.1418$ &   $0.0181$ &   $0.0099$ \\
		\hline
		Parameters of the  & & &  \\
		droop function $K_{i,k}$,  & & & \\
		$i\in\mathbb{N}_{[1,5]}$ in \eqref{eq:droop_K}  & & & \\
		according to $d_k$: & & & \\
		$[w_{\max,i,k}]$ & $\left[\begin{array}{c}29613  \\  45830  \\  33138 \\   37369 \\   25383\end{array}\right]$ & $\left[\begin{array}{c} 232320  \\  359540  \\  259980 \\   293160  \\  199130\end{array}\right]$& $\left[ \begin{array}{c} 423930  \\  656080  \\  474390 \\   534950  \\  363360\end{array} \right]$  \\
		$w_{\min,i,k}=-w_{\max,i,k}$ & & & \\
		$w_{m,i,k}=w_{n,i,k}$ & $0$ & $0$ & $0$ \\
 		\hline 
	\end{tabular} 
\end{table}

\begin{figure}[h!] 
	\begin{center}
		\includegraphics[scale=0.8]{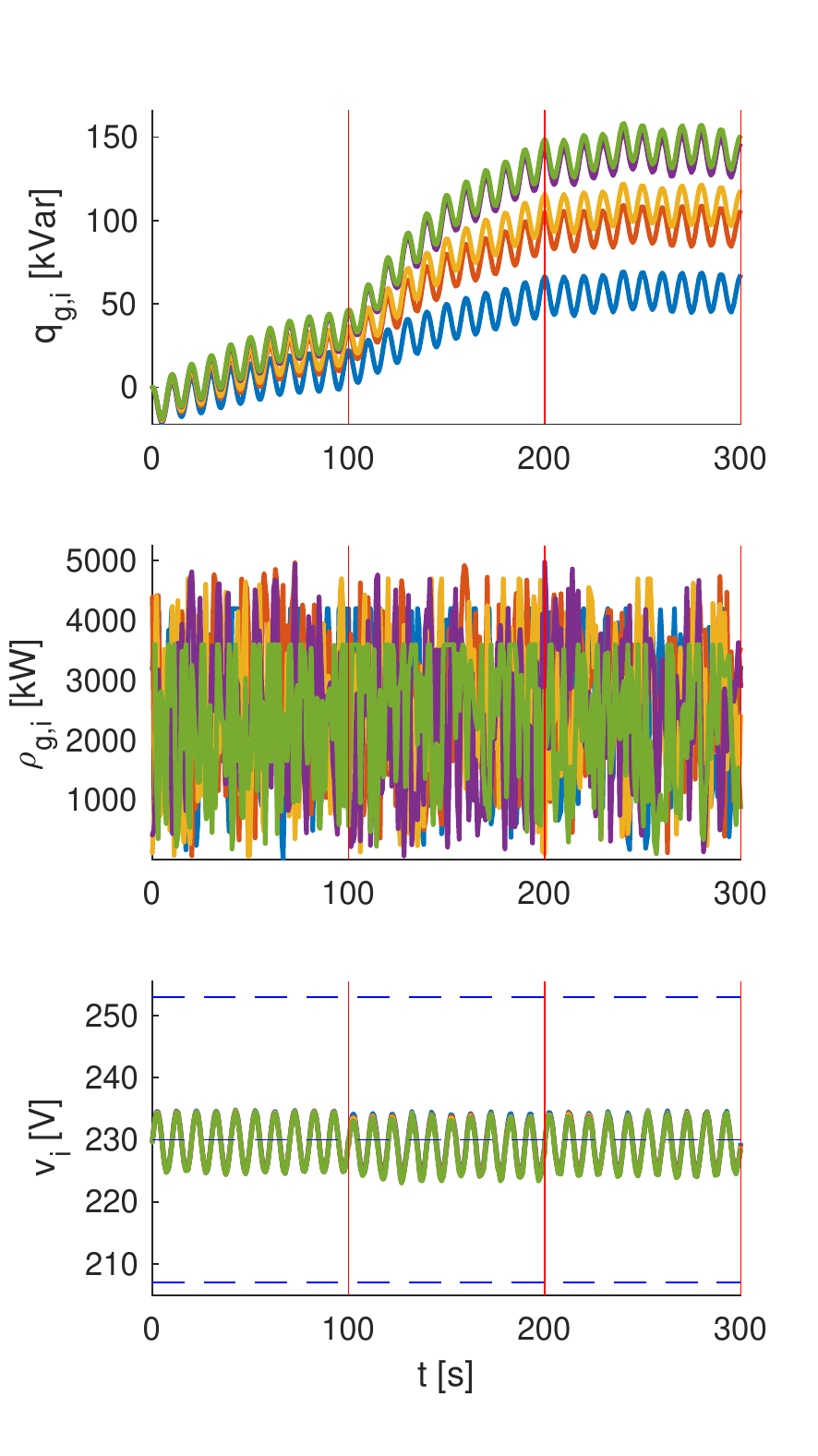} \vspace{-2em}
		\caption{Non-adaptive droop controllers: Vertical {\color{red}red} lines indicate the observation instants $t_k$, $k\in\mathbb{N}_{[0,2]}$. Legend: Customer $1$ (blue), $2$ (red), $3$ (orange), $4$ (purple), $5$ (green).  (Top) Solid lines indicate the output of the droop controllers $q_{g,i}$. (Middle) The active power $\rho_{g,i}$ injected into the distribution line by each customer $i$. (Bottom) Solid lines indicate the voltages of each customer $v_i$ are within $10$ \% of the nominal voltage $\bar{v}$ for all $t\geq 0$. Dashed lines indicate $(1\pm p)\bar{v}$, where $p=0.1$. } \label{fig:non_time}
	\end{center}
\end{figure}

For comparison, we consider the case where we do not adapt the droop controller \eqref{eq:droop} in accordance to the known bounds $\Delta_{c,k}$ and $\Delta_{\rho,k}$ for each observation window. Instead, only an initial design based on the largest known bound for the entire observation window $[t_0,t_3)$ for the droop controller \eqref{eq:droop} is used, i.e. we design the droop function $K_{i,k}$ for all $i\in\mathbb{N}_{[1,N]}$, $k\in\mathbb{N}_{[0,3]}$ based on the maximal known bounds $\Delta_{c}=\underset{k\in\mathbb{N}_{[0,2]}}{\max} \Delta_{c,k} = 1200$, $\Delta_{\rho}=\underset{k\in\mathbb{N}_{[0,2]}}{\max} \Delta_{\rho,k} = 1300$, resulting in $d_{i,k}=0.3251$ and consequently,  $w_{m,i,k}=w_{n,i,k}=0$, $\left[w_{\max,i,k}\right]=\left[12918,    19992,    14456,    16301,    11073
 \right]$, and $w_{\min,i,k}=-w_{\max,i,k}$ for $i\in\mathbb{N}_{[1,5]}$ and $k\in\mathbb{N}_{[0,3]}$. Figure \ref{fig:non_time} shows that the voltage regulation objective \eqref{eq:restate_voltage_bounds} is achieved.

We use the performance metrics $\mathcal{P}_{q}:=\int_{0}^{t_3} |q_{g}(s)| ds$ and $\mathcal{P}_{\rho}:=\int_{0}^{t_3} |\rho_{g}(s)| ds$ to compare the adaptive vs. non-adaptive droop controllers. The metric $\mathcal{P}_{q}$ is motivated by achieving the voltage regulation objective with less reactive power generated by the inverter.  The metric $\mathcal{P}_{\rho}$ is used to monitor profit as increasing the injection of active power into the distribution network increases the generation of revenue. For the simulation scenario described, the adaptive scheme outperforms the non-adaptive scheme in the aforementioned performance metrics with $\mathcal{P}_{q}=2.9456$ [kJ], $\mathcal{P}_{\rho}=1.6182$ [MJ] (adaptive) and $\mathcal{P}_{q}=30.833$ [kJ], $\mathcal{P}_{\rho} = 1.2404$ [MJ] (non-adaptive).



\section{Conclusions and future work} \label{sec:conc}
We have presented a scheme that adapts to the known bounds of the load's net power consumption with guarantees on each customers' voltage being within a desired safety margin through local inverters injecting reactive power into the distribution network. Future work will take into account the unmodelled dynamics in the power distribution line.

\section{Acknowledgement}
We are thankful to Catalin Gavriluta and David Umsonst for the discussions leading to this work and for reproducing Figure \ref{fig:cigre} from \cite{{strunz2009benchmark}}. We are also grateful to the constructive comments given by the anonymous reviewers.  

\bibliography{grid_model_analysis.bib}             

\appendix


\section{Proof of Theorem \ref{thm:voltage_reg}} 
Our aim is to show that given $\epsilon \geq 0$, $\{\Delta_{\rho,k}\}_{k\in\mathbb{N}}$ and $\{\Delta_{c,k}\}_{k\in\mathbb{N}}$, $|y_i(t)|\leq \epsilon$, for all $i\in\mathbb{N}_{[1,N]}$ and for all $t\geq 0$. Since 
\begin{equation}
	\begin{array}{lll}
		\left| y_i(t) \right| & \leq & \underset{i\in\mathbb{N}_{[1,N]}}{\max}\left|y_i(t) \right| =: \left| y(t) \right| \\
		& \leq & \left| H q_{g} \right| + \left| \phi(\rho,q_c) \right| + |y_0|, 
	\end{array} \label{eq:y_inequality}
\end{equation}
we derive a bound for each term as follows.

\begin{enumerate}[(i)]
	\item \textbf{Bound on $\left| \phi(\rho,q_c) \right|$}: \\
			\begin{equation} \label{eq:bound_phi}
				\begin{array}{lll} 
						\left| \phi(\rho,q_c) \right| & \leq & \underset{i\in\mathbb{N}_{[0,N-1]}}{\sum} \left| \phi_i(\rho,q_c) \right| \\
						&& + \underset{i\in\mathbb{N}_{[0,N-2]}}{\sum} 2 \left| \beta'_{i}(\rho_{i+1},q_{c,i+1}) \right|.
				\end{array}
			\end{equation}
			For $t\in[t_{k},t_{k+1})$, $k\in\mathbb{N}$, we have 
			\begin{equation}
				\begin{aligned}
						\underset{i\in\mathbb{N}_{[0,N-1]}}{\sum} &\left| \phi_i(\rho,q_c) \right|  \\ &\leq \underset{i\in\mathbb{N}_{[0,N-1]}}{\sum} 2 R_i \left( \underset{j\in\mathbb{N}_{[i+1,N]}}{\sum} \left| \rho_j \right| \right)\\
						& \; + \underset{i\in\mathbb{N}_{[0,N-1]}}{\sum} 2 X_i \left( \underset{j\in\mathbb{N}_{[i+1,N]}}{\sum} \left| q_{c,j} \right| \right) \\
						& \; + \underset{i\in\mathbb{N}_{[0,N-1]}}{\sum} 2 R'_i \left| \rho_{i+1} \right| + 2 X'_{i} \left| q_{c,i+1} \right| \\ 
						& \leq  \underset{i\in\mathbb{N}_{[0,N-1]}}{\sum} 2 \bar{R} \left( \underset{j\in\mathbb{N}_{[i+1,N]}}{\sum} \Delta_{\rho,k} \right) \\
						& \; + \underset{i\in\mathbb{N}_{[0,N-1]}}{\sum} 2 \bar{X} \left( \underset{j\in\mathbb{N}_{[i+1,N]}}{\sum} \Delta_{q,k} \right) \\
						&\; +  \underset{i\in\mathbb{N}_{[0,N-1]}}{\sum} 2 \bar{R}' \Delta_{\rho,k} + 2 \bar{X}' \Delta_{c,k} \\
						& =  2 \underset{i\in\mathbb{N}_{[0,N-1]}}{\sum} \left( N - i \right) \left( \bar{R} \Delta_{\rho,k} + \bar{X} \Delta_{q,k} \right) \\
						& \;+ 2N \left( \bar{R}' \Delta_{\rho,k} + \bar{X}' \Delta_{q,k} \right)\\
						& =  N(N+1) \left( \bar{R} \Delta_{\rho,k} + \bar{X} \Delta_{q,k} \right) \\
						& \; + 2N \left( \bar{R}' \Delta_{\rho,k} + \bar{X}' \Delta_{q,k} \right),
				\end{aligned} \label{eq:bound_phi_1}
			\end{equation}
			and 
			\begin{equation}
	\begin{aligned}
			\underset{i\in\mathbb{N}_{[0,N-2]}}{\sum} & \beta'_{i}\left(\rho_{i+1}, q_{c,i+1} \right) \\
			 & \leq  \underset{i\in\mathbb{N}_{[0,N-2]}}{\sum} 2 R'_{i} \left|\rho_{i+1}\right| + 2 X'_{i} \left|q_{c,i+1}\right| \\
			& \leq  \underset{i\in\mathbb{N}_{[0,N-2]}}{\sum} 2 \bar{R}' \Delta_{\rho,k} + 2 \bar{X}' \Delta_{q,k} \\
			& =  2 (N-1) \left(\bar{R}' \Delta_{\rho,k} + \bar{X}' \Delta_{q,k} \right), 	
	\end{aligned} \label{eq:bound_beta}
\end{equation}
where $\bar{R}:= \max_{i} R_i$, $\bar{X}:= \max_i X_i$, $\left|\rho_i(t)\right| \leq \Delta_{\rho,k}$, $\left|q_{c,i}(t)\right| \leq \Delta_{q,k}$ for $t\in[t_{k},t_{k+1})$ and we obtain the last equality of \eqref{eq:bound_phi_1} due to 
${\underset{i\in\mathbb{N}_{[0,N-1]}}{\sum} N-i = \frac{1}{2}N(N+1)}$.

Therefore, we obtain the following from \eqref{eq:bound_phi}, \eqref{eq:bound_phi_1} and \eqref{eq:bound_beta}.
\begin{equation}
	\begin{array}{lll}
		 \left| \phi(\rho,q_c) \right| & \leq & N(N+1) \left( \bar{R} \Delta_{\rho,k} + \bar{X} \Delta_{q,k} \right) \\
		 								&& + 2(2N-1) \left( \bar{R}' \Delta_{\rho,k} + \bar{X}' \Delta_{q,k} \right) =: \Delta_{\phi,k}.
	\end{array} \label{eq:bound_phi_final}
\end{equation} 
			
	\item \textbf{Bound on $\left| H q_{g} \right|$}: \\
			For $t\in[t_{k},t_{k+1})$, $k\in\mathbb{N}$, by noting that $q_{g}(t_k)=0$, the solution to \eqref{eq:model_compact} is
			\begin{equation}
				\begin{array}{lll}
					Hq_{g}(t) & = & - \int_{t_{k}}^{t} \exp(A(\tau)(t-s)) H A(\tau) K_{k}(y(s)) \, ds. 
				\end{array}
			\end{equation}
			Since $A(\tau)$ is Hurwitz, we have the following properties:
			\begin{itemize}
					\item $\left| \exp\left( A(\tau) (t-t_k-s) \right) \right| \leq \exp \left(-\frac{1}{\tau_{\max}} (t-t_k-s) \right)$,
					\item $\left|-A(\tau)\right| \leq \frac{1}{\tau_{\min}}$.
			\end{itemize}
			Further, by the slope-restriction property of $K_k(y)$ as stated in Assumption \ref{assum:sector_bound}, we have $\left|K_{k}(y) \right|\leq d_k \left| y \right| \leq d_k \left( \left|H q_{g}\right| + \left| \phi(\rho,q_c) \right| + \left| y_0 \right|  \right)$. Therefore,
			\begin{equation}
				\begin{array}{lll}
				\left| H q_g(t) \right| & \leq & \frac{1}{\tau_{\min}} d_k \left| H \right| \\
				&& \times \int_{t_k}^{t} \exp\left(-\frac{1}{\tau_{\max}} (t-t_k-s) \right) \\
				&& \qquad \times \left( \left|H q_{g}\right| + \left| \phi(\rho,q_c) \right| + \left| y_0 \right|  \right)  ds.
				\end{array}
			\end{equation}
			From \eqref{eq:bound_phi_final} and $\left| y_0 \right|:= \left| \bar{v}^2 - v'^{2}_{0} \right| \leq \epsilon_y$, we obtain
			\begin{equation}
				\begin{array}{lll}
					\left| H q_g(t) \right| & \leq & d_k \frac{\tau_{\max}}{\tau_{\min}} \left|H\right| \left( \Delta_{\phi,k}  + \epsilon_{y} \right)  \\
					&& \times \int_{t_k}^{t} \exp\left(- \frac{1}{\tau_{\max}} \left( t-t_k- s \right) \right) \\
					&& + \frac{d_k}{\tau_{\min}} \left| H \right|  \int_{t_k}^{t} \exp\left(- \frac{1}{\tau_{\max}} \left( t-t_k- s \right) \right)\\
					&& \qquad \qquad \qquad \qquad \times  \left| H q_{g}(s) \right| ds	
				\end{array}
			\end{equation}
			We employ the boundedness property of $K_{k}(y)$ (i.e. $\left|K_{k}(y)\right| \leq \bar{K}_{k}$) to obtain
			\begin{equation}
				\begin{array}{lll}
					\left| H q_g(t) \right|& \leq &  d_k \left(\frac{\tau_{\max}}{\tau_{\min}} \left|H\right| \left( \Delta_{\phi,k}  + \epsilon_{y} \right) \right. \\
					&& \left. \qquad \qquad + \left(\frac{\tau_{\max}}{\tau_{\min}} \right)^2 \left| H \right|^2 \bar{K}_{k} \right) \\
					&& \times \int_{t_k}^{t} \exp\left(- \frac{1}{\tau_{\max}} \left( t-t_k- s \right) \right) ds \\
					& \leq & d_k \left( \frac{\tau_{\max}}{\tau_{\min}} \left|H\right| \left( \Delta_{\phi,k}  + \epsilon_{y} \right)  \right. \\
					&&  \qquad \left. + \left(\frac{\tau_{\max}}{\tau_{\min}} \right)^2 \left| H \right|^2 \bar{K}_{k} \right) \exp\left( \frac{t_{k}}{\tau_{\max}}  \right) \\
					
				\end{array} \label{eq:hq_bound}
			\end{equation}
			where we obtain the last inequality due to 
			\begin{equation}
			\begin{aligned} 
				\int_{t_k}^{t} & \exp\left(\frac{1}{\tau_{\max}}\left(t-t_k-s\right) \right) ds \\
				& = \tau_{\max}\exp\left( t_k / \tau_{\max} \right) \left(1-\exp\left( -\left(t- t_k \right) / \tau_{\max} \right) \right) \\
				& \leq \tau_{\max}\exp\left( t_k / \tau_{\max} \right).
				\end{aligned}
			\end{equation}
			\item \textbf{Conclusion}: From \eqref{eq:y_inequality}, \eqref{eq:bound_phi_final}, \eqref{eq:hq_bound} and \eqref{eq:cond_volt_reg}, we have shown that $\left| y_i(t) \right| \leq \epsilon$ for $i\in\mathbb{N}_{[1,N]}$ and for all $t\geq 0$ as desired.
\end{enumerate}
\hfill $\Box$

\end{document}